\documentclass[fleqn]{mat01}
\usepackage{times,mathtimy,amssymb,latexsym}
\begin{document}

\setcounter{page}{325} \firstpage{325}

\newtheorem{theore}{Theorem}
\renewcommand\thetheore{\arabic{section}.\arabic{theore}}
\newtheorem{theor}[theore]{\bf Theorem}
\newtheorem{lem}[theore]{Lemma}
\newtheorem{coro}[theore]{\rm COROLLARY}
\newtheorem{propo}[theore]{\rm PROPOSITION}

\def\exam{\trivlist\item[\hskip\labelsep{\it Example.}]}

\title{Inequalities involving upper bounds for certain matrix
operators}

\markboth{R~Lashkaripour and D~Foroutannia}{Upper bounds for
certain matrix operators}

\author{R~LASHKARIPOUR and D~FOROUTANNIA}

\address{Department of Mathematics, Sistan and Baluchestan
University, Zahedan, Iran\\
\noindent E-mail: lashkari@hamoon.usb.ac.ir;
d$_{-}$foroutan@math.com}

\volume{116}

\mon{August}

\parts{3}

\pubyear{2006}

\Date{MS received 15 February 2006}

\begin{abstract}
In this paper, we considered the problem of finding the upper
bound  Hausdorff matrix operator from sequence spaces $l_p(v)$ (or
$d(v,p)$) into $l_p(w)$ (or $d(w,p)$). Also we considered the
upper bound problem for matrix operators from $d(v,1)$ into
$d(w,1)$, and matrix operators from $e(w,\infty)$ into
$e(v,\infty)$, and deduce upper bound for Cesaro, Copson and
Hilbert matrix operators, which are recently considered in
\cite{5} and \cite{6} and similar to that in \cite{10}.
\end{abstract}

\keyword{Inequality; norm; summability matrix; Hausdorff matrix;
Hilbert matrix; weighted sequence space; Lorentz sequence space.}

\maketitle

\section{Introduction}

We study the norm of a certain matrix operator on $l_p(w)$ and
Lorentz sequence spaces $d(w,p)$, $p\geq 1$, which is considered
in \cite{2} on $l_p$ spaces and in \cite{6,7,8} and \cite{9} on
$l_p(w)$ and $d(w,p)$ for some matrix operator such as Cesaro,
Copson and Hilbert operators.

Let $l_p$ be the normed linear space of all sequences $x=(x_n)$
with finite norm $\|x\|_p$, where
\begin{equation*}
\|x\|_p = \left( \sum_{n=1}^{\infty}|x_n|^p \right)^{1/p}.
\end{equation*}
Let $w=(w_n)$ be a sequence with positive entries. For $p\ge 1$,
we define the weighted sequence space $l_p(w)$ as follows:
\begin{equation*}
l_p(w) = \left\{(x_n)\hbox{:}\ \sum_{n=1}^{\infty}
w_n|x_n|^p<\infty\right\},
\end{equation*}
with norm, $\|\cdot\|_{p,w}$, which is defined as follows:
\begin{equation*}
\|x\|_{w,p}= \left( \sum_{n=1}^{\infty}w_n|x_n|^p \right)^{1/p}.
\end{equation*}
Also, if $w=(w_n)$ is a decreasing sequence of positive number
such that $\lim_{n\rightarrow\infty} w_n=0$ and
$\sum_{n=1}^{\infty}w_n=\infty$, then the Lorentz sequence space
$d(w,p)$ is defined as follows:
\begin{equation*}
d(w,p)= \left\{ (x_n)\hbox{:}\ \sum_{n=1}^{\infty}
w_nx^{*p}_{n}<\infty\right\},
\end{equation*}
where $(x^*_n)$ is the decreasing rearrangement of $(|x_n|)$. In
fact $d(w,p)$ is the space of null sequences $x$ for which $x^*$
is in $l_p(w)$, with norm $\|x\|_{d(w,p)}=\|x^*\|_{w,p}$.

Let $X^*_k=x^*_1+\cdots+x^*_k$ and $W_k=w_1+\cdots+w_k$. We define
the weighted sequence space $e(w,\infty)$ as follows:
\begin{equation*}
e(w,\infty) = \left\{(x_n)\hbox{:}\ \sup_k \frac{X^*_k}{W_k}
<\infty \right\},
\end{equation*}
with norm $\|\cdot\|_{w,\infty}$, which is defined as follows:
\begin{equation*}
\|x\|_{w,\infty}=\sup_k\frac{X^*_k}{W_k}.
\end{equation*}
Our objective in \S2 is to give a generalization of some results
obtained by Bennett \cite{1,2} and Jameson and Lashkaripour
\cite{6}, for Hausdorff matrix operators on the weighted sequence
space. In \S3 we try to solve the problem of finding the norm of
matrix operators from $d(v,1)$ into $d(w,1)$, and matrix operators
from $e(w,\infty)$ into $e(v,\infty)$, and we deduce upper bound
for certain matrix operators such as Cesaro, Copson and Hilbert
operators.

\section{Hausdorff matrix operator on $\pmb{l_p(w)}$ and
$\pmb{d(w,p)}$}

In this section, we consider the Hausdorff matrix  operator
$H(\mu)=(h_{j,k})$, such that
\begin{align*}
h_{j,k} = \left\{\begin{array}{l@{\quad}l} \left(
\begin{smallmatrix}
j-1\\[.1pc]
k-1
\end{smallmatrix}\right)
\triangle^{j-k}a_k, &\mbox{if}\ \ 1\le k\leq j,\\[.4pc]
0, &\mbox{if}\ \ k>j,
\end{array} \right.
\end{align*}
where $\triangle$ is the difference operator; that is
\begin{equation*}
\triangle a_k=a_k-a_{k+1},
\end{equation*}
and $(a_k)$ is a sequence of real numbers, normalized so that
$a_1=1$.

If
\begin{equation*}
a_k = \int_0^1\theta^{k-1}\hbox{\rm d}\mu(\theta),\quad
k=1,2,\dots,
\end{equation*}
where $\mu$ is a probability measure on $[0,1]$, then for all $j,
k=1,2,\dots$, we have
\begin{equation*}
h_{j,k} = \left\{ \begin{array}{l@{\quad}l} \left(\begin{smallmatrix}
j-1\\[.1pc]
k-1
\end{smallmatrix} \right)
\int_0^1 \theta^{k-1}(1-\theta)^{j-k}\hbox{\rm d}\mu(\theta),
&\mbox{if}\ \ 1\le k\leq j\\[.4pc]
0, &\mbox{if}\ \ k>j
\end{array}.\right.
\end{equation*}
The Hausdorff matrix contains the famous classes of matrices.
These classes are as\break follows:
\begin{enumerate}
\renewcommand\labelenumi{(\roman{enumi})}
\leftskip .4pc
\item Choice $\hbox{\rm d}\mu(\theta)=\alpha(1-\theta)^{\alpha-1}
\hbox{\rm d}\theta$ gives the Cesaro matrix of order $\alpha$.

\item Choice $\hbox{\rm d}\mu(\theta)= \hbox{point evaluation at} \
\theta=\alpha$ gives the Euler matrix of order $\alpha$.

\item Choice $\hbox{\rm d}\mu(\theta)=
\frac{|\log\theta|^{\alpha-1}}{\Gamma(\alpha)} \hbox{\rm d}\theta$
gives the H\"{o}lder matrix of order $\alpha$.

\item Choice $\hbox{\rm d}\mu(\theta)=\alpha\theta^{\alpha-1}
\hbox{\rm d}\theta$ gives the Gamma matrix of order $\alpha$.
\end{enumerate}
The Cesaro, H\"{o}lder and Gamma matrices have non-negative
entries, whenever $\alpha>0$. Also the Euler matrix is
non-negative, when $0\leq\alpha\leq 1.$ So, if we obtain the norm
of the Hausdorff matrix, then it is also an upper bound for the
above matrices.

Now consider the operator $A$ defined by $Ax=y$, where
$y_i=\sum_{i =1}^{\infty}a_{i,j}x_j$. We write $\|A\|_{v,w,p}$ for
the norm of $A$ as an operator from $l_p(v)$ into $l_p(w)$, and
$\|A\|_{w,p}$ for the norm of $A$ as an operator from  $l_p(w)$
into itself, and $\|A\|_{p}$ for the norm of $A$ as an operator
from $l_p$ into itself, and  $\|A\|_{d(w,p)}$ for the norm of $A$
as an operator from ${d(w,p)}$ into itself.

The following conditions are needed to convert statements for
$l_p(w)$ to ones for ${d(w,p)}$. We assume throughout that
\begin{enumerate}
\renewcommand\labelenumi{(\arabic{enumi})}
\leftskip .15pc
\item For all $i, j, a_{i,j}\geq0$.

\item For all subsets $M, N$ of natural numbers having  $m, n$
elements respectively, we have
\begin{equation*}
\hskip -1.25pc \sum_{i\in M}\sum_{j\in N}a_{i,j}\le\sum_{i=1}^m
\sum_{j=1}^na_{i,j}.
\end{equation*}

\item $\sum_{i=1}^{\infty}w_ia_{i,1}$ is convergent.

Condition (1) implies that $|A(x)|\leq A(|x|)$  and hence the
non-negative sequences are sufficient to determine norm of $A$.
\end{enumerate}

\begin{propo}$\left.\right.$\vspace{.5pc}

\noindent Let $p\ge 1$ and $A=(a_{i,j})$ be an operator with
conditions $(1)$ and $(2)$. Then
\begin{equation*}
\|A(x)\|_{d(w,p)}\le\|A(x^*)\|_{d(w,p)},
\end{equation*}
for all non-negative elements $x$ of $d(w,p)$. Hence
decreasing{\rm ,} non-negative elements are sufficient to
determine norm of $A$.
\end{propo}

Condition (3) ensured that at least finite sequence are mapped
into $d(w,1)$.

\begin{propo} {\rm (Lemma~1 of \cite{5})}$\left.\right.$\vspace{.5pc}

\noindent Let $p\ge 1$ and $A=(a_{i,j})$ be an operator with
non-negative entries. Also{\rm ,} let $A$ map $d(w,p)$ into
itself. If for $x\in d(w,p)${\rm ,} we set $Ax=y$ such that
$y_i=\sum_{j=1}^{\infty}a_{i,j}x_j${\rm ,} then the following
conditions are equivalent{\rm :}

\begin{enumerate}
\renewcommand\labelenumi{\rm (\alph{enumi})}
\leftskip .15pc
\item $y_1\ge y_2\ge\cdots\ge 0$ when $x_1\ge x_2\ge\cdots\ge 0$.

\item $r_{i,n}=\sum_{j=1}^{n}a_{i,j}$ decreases with $i$ for each
$n$.
\end{enumerate}
\end{propo}

In the following statement, we assume $(v_n)$ and $(w_n)$ to be
non-negative decreasing sequences with $v_1=1$.

\setcounter{theore}{0}
\begin{theor}[\!]
Let $H(\mu)$ be the Hausdorff matrix operator and $p>1$. Then the
Hausdorff matrix operator maps $l_p(v)$ into $l_p(w)${\rm ,} and
\begin{align*}
&\left( \inf{\frac{w_n}{v_n}} \right)^{1/p} \int_0^1\theta^{-1/p}
\hbox{\rm d}\mu(\theta)\\[.3pc]
&\quad\, \leq \|H\|_{v,w,p}\leq \left( \sup{\frac{w_n}{v_n}}
\right)^{1/p}\int_0^1 \theta^{-1/p}\hbox{\rm d}\mu(\theta).
\end{align*}
Therefore the Hausdorff matrix operator maps $l_p(w)$ into
itself{\rm ,} and
\begin{equation*}
\|H\|_{w,p}=\int_0^1\theta^{-1/p}\hbox{\rm d}\mu(\theta).
\end{equation*}
\end{theor}

\begin{proof}
Let $x$ be a non-negative sequence. Since $(w_n)$ is decreasing,
and applying Theorem~216 of \cite{3}, we have
\begin{align*}
\|Hx\|_{w,p}^p &= \sum_{j=1}^{\infty}w_j \left(\sum_{k=1}^j
\left(\begin{array}{l} j-1\\[.1pc] k-1
\end{array}\right)
\left(\int_0^1
\theta^{k-1}(1-\theta)^{j-k}\hbox{\rm d}\mu(\theta)\right)x_k\right)^p\\[.4pc]
&\leq \sum_{j=1}^{\infty} \left(\sum_{k=1}^j \left(\begin{array}{l} j-1\\
k-1
\end{array}\right)
\left(\int_0^1 \theta^{k-1}(1-\theta)^{j-k}\hbox{\rm d}\mu(\theta)\right)w_k^{1/p}x_k\right)^p\\[.4pc]
&\leq \left(\int_0^1\theta^{-1/p}\hbox{\rm d}\mu(\theta)\right)^{p}\sum_{j=1}^{\infty}w_jx_j^p\\[.4pc]
&= \left( \int_0^1\theta^{-1/p}\hbox{\rm d}\mu(\theta)\right)^{p}\sum_{j=1}^{\infty}\frac{w_j}{v_j}v_jx_j^p\\[.4pc]
&\leq \sup\frac{w_j}{v_j} \left(\int_0^1\theta^{-1/p}\hbox{\rm d}
\mu(\theta) \right)^p\|x\|_{v,p}^p.
\end{align*}
Hence
\begin{equation*}
\|Hx\|_{w,p}\le \left(\sup \frac{w_n}{v_n} \right)^{1/p}
\left(\int_0^1\theta^{-1/p}\hbox{\rm
d}\mu(\theta)\right)\|x\|_{v,p},
\end{equation*}
and so
\begin{equation*}
\|H\|_{v,w,p}\le \left( \sup\frac{w_n}{v_n}
\right)^{1/p}\int_0^1\theta^{-1/p}\hbox{\rm d}\mu(\theta).
\end{equation*}
It remains to prove the left-hand inequality. We take
\begin{equation*}
0<\delta<\frac{1}{p},\quad x_n=(n)^{-\frac{1}{p}-\delta}
\end{equation*}
and any positive $\epsilon$, where $0<\epsilon< 1$; and choose
$\alpha$, $N$, and $\delta$ such that
\begin{align*}
&\left(1+\frac{1}{\alpha}\right)^{-2/p} > 1-\epsilon,\\[.4pc]
&\int_{\alpha/n}^1\theta^{-1/p}\hbox{\rm d}\mu(\theta)
> (1-\epsilon)\int_0^1\theta^{-1/p}\hbox{\rm d}\mu(\theta), \quad n\ge N,\\[.4pc]
&\sum_{n=N}^{\infty}w_nx_n^p >
(1-\epsilon)\sum_{n=1}^{\infty}w_nx_n^p.
\end{align*}
Since $(x_n)\in l_p$, and $0<v_n\le 1$, we deduce that $(x_n)\in
l_p(v)$. Also, we have
\begin{align*}
(Hx)_n &= \sum_{m=1}^{n}\left(\begin{array}{l} n-1\\[.1pc] m-1
\end{array}\right)
\left( \int_0^1
\theta^{m-1}(1-\theta)^{n-m}\hbox{\rm d}\mu(\theta)\right) x_m\\[.4pc]
&\ge(1-\epsilon)^2x_n\int_0^1\theta^{-1/p}\hbox{\rm d}\mu(\theta),
\quad n\ge N,
\end{align*}
and so
\begin{equation*}
w_n^{1/p}(Hx)_n\ge(1-\epsilon)^2w_n^{1/p}x_n\int_0^1\theta^{-1/p}
\hbox{\rm d}\mu(\theta), \quad n\ge N.
\end{equation*}
Hence
\begin{align*}
\|Hx\|_{w,p}^p &\ge \sum_{n=N}^{\infty}w_n(Hx)_n^p\\[.4pc]
&\ge (1-\epsilon)^{2p}\left( \int_0^1\theta^{-1/p}\hbox{\rm
d}\mu(\theta)\right)^p
\sum_{n=N}^{\infty}w_nx_n^p\\[.4pc]
&\ge (1-\epsilon)^{2p+1}\left(\int_0^1\theta^{-1/p}
\hbox{\rm d}\mu(\theta)\right)^p \sum_{n=1}^{\infty}w_nx_n^p\\[.4pc]
&= (1-\epsilon)^{2p+1}\left( \int_0^1\theta^{-1/p}\hbox{\rm
d}\mu(\theta) \right)^p
\sum_{n=1}^{\infty}\frac{w_n}{v_n}v_nx_n^p\\[.4pc]
&\ge \inf\frac{w_n}{v_n}(1-\epsilon)^{2p+1}\left(
\int_0^1\theta^{-1/p}\hbox{\rm d} \mu(\theta)\right)^p
\|x\|_{v,p}^p.
\end{align*}
Since $\epsilon$ is arbitrary, if $\epsilon\longrightarrow 0$, we
have
\begin{equation*}
\|Hx\|_{w,p}^p\ge\inf\frac{w_n}{v_n} \left(
\int_0^1\theta^{-1/p}\hbox{\rm d}\mu(\theta)
\right)^p\|x\|_{p,v}^p,
\end{equation*}
and this completes the proof of the theorem.\hfill $\Box$
\end{proof}

\setcounter{theore}{0}
\begin{coro}$\left.\right.$\vspace{.5pc}

\noindent Let $p>1$ and $p^*=\frac{p}{p-1}$. Then Cesaro{\rm ,}
H\"{o}lder{\rm ,} Gamma and Euler operators map $l_p(w)$ into
$l_p(w)$. Also{\rm ,} we have
\begin{align*}
\|C(\alpha)\|_{w,p}&=\frac{\Gamma(\alpha+1)\Gamma(1/p^*)}{\Gamma\left(\alpha+\frac{1}{p^*}\right)},\quad
\alpha>0;\\[.4pc]
\|H(\alpha)\|_{w,p}&=\frac{1}{\Gamma(\alpha)}\int_0^1\theta^{-\frac{1}{p}}|\log\theta|^{\alpha-1}\hbox{\rm d}\theta,\quad \alpha>0;\\[.4pc]
\|\Gamma(\alpha)\|_{w,p}&=\frac{\alpha p}{\alpha p-1},\quad \alpha p>1;\\[.4pc]
\|E(\alpha)\|_{w,p}&=\alpha^{-1/p},\quad 0<\alpha<1.
\end{align*}
\end{coro}

\begin{proof}
It is elementary.\hfill $\Box$
\end{proof}

The following corollary is an extension of Theorem~326 (p.~239 of
\cite{4}).

\begin{coro}$\left.\right.$\vspace{.5pc}

\noindent If $x$ and $w$ are non-negative sequences and $w$ is
decreasing{\rm ,} then
\begin{equation*}
\sum_{n=1}^{\infty}w_n \left(\frac{1}{n}\sum_{i=1}^{n}x_i
\right)^p\le {p^*}^p\left(\sum_{n=1}^{\infty}w_nx_n^p\right).
\end{equation*}
\end{coro}

\begin{proof}
For Cesaro operator we apply Corollary~2.1 with $\alpha=1$.\hfill
$\Box$
\end{proof}

\begin{coro}$\left.\right.$\vspace{.5pc}

\noindent If $H(\mu)$ is the Hausdorff matrix operator on $l_p$
and $p>1$, then
\begin{equation*}
\|H\|_{p}=\int_0^1\theta^{-1/p}\hbox{\rm d}\mu(\theta).
\end{equation*}
\end{coro}

\begin{proof}
By taking $w_n=1$ for all $n$, we have the corollary.\hfill $\Box$
\end{proof}

\begin{coro}$\left.\right.$\vspace{.5pc}

\noindent Let $p>1$. Then Cesaro{\rm ,} H\"{o}lder{\rm ,} Gamma
and Euler operators map $l_p$ into $l_p$. Also{\rm ,} we have
\begin{align*}
\|C(\alpha)\|_{p}&=\frac{\Gamma(\alpha+1)\Gamma(1/p^*)}{\Gamma\left(\alpha+\frac{1}{p^*}\right)},\quad \alpha>0;\\[.4pc]
\|H(\alpha)\|_{p}&=\frac{1}{\Gamma(\alpha)}\int_0^1\theta^{-\frac{1}{p}}|\log\theta|^{\alpha-1}\hbox{\rm d}\theta,\quad \alpha>0;\\[.4pc]
\|\Gamma(\alpha)\|_{p}&=\frac{\alpha p}{\alpha p-1},\quad \alpha p>1;\\[.4pc]
\|E(\alpha)\|_{p}&=\alpha^{-1/p},\quad 0<\alpha<1.
\end{align*}
\end{coro}

\begin{proof}
It is  elementary.\hfill $\Box$
\end{proof}

\setcounter{theore}{1}
\begin{theor}[\!]
Let $p>1$ and $H(\mu)$ be the Hausdorff matrix operator with
condition $(2)$. Then the Hausdorff matrix operator{\rm ,}
$H(\mu)${\rm ,} maps $d(w,p)$ into itself{\rm ,} and we have
\begin{equation*}
\|H\|_{d(w,p)}=\int_0^1\theta^{-1/p}{\rm d}\mu(\theta).
\end{equation*}
\end{theor}

\begin{proof}
By Propositions~2.1 and 2.2, it is enough to consider non-negative
decreasing sequences. For such sequences, we have
$\|Hx\|_{d(w,p)}=\|Hx\|_{w,p}$, and so applying Theorem~1.1, we
deduce the theorem.\hfill $\Box$
\end{proof}

\begin{exam}{\rm
Suppose $p>1$. Since $\Gamma(1)=C(1)$ and they satisfy condition
(2), we have
\begin{equation*}
\|\Gamma(1)\|_{d(w,p)}=\|C(1)\|_{d(w,p)}=p^*.
\end{equation*}
Also
\begin{align*}
C(2)=\left[\begin{array}{l@{\quad}l@{\quad}l@{\quad}l@{\quad}l@{\quad}l}
1 &0\\[.1pc]
2/3 &1/3 &0\\[.1pc]
3/6 &2/6 &1/6 &0\\[.1pc]
4/10 &3/10 &2/10 &1/10 &0\\[.1pc]
\cdot &\cdot &\cdot &\cdot &\cdot &\cdot
\end{array}\right]
\end{align*}
has condition (2) and so $\|C(2)\|_{d(w,p)}=p^*(2p)^*$.}
\end{exam}

\section{Matrix operators on $\pmb{d(w,1})$ and
$\pmb{e(w,\infty)}$}

Here we consider the upper bound problem for matrix operators from
$d(v,1)$ into $d(w,1)$, and matrix operators from $e(w,\infty)$
into $e(v,\infty)$. If $x\in d(w,1)$, we denote norm of $x$ with
$\|x\|_{w,1}$  and if $x\in e(w,\infty)$, we denote norm of $x$
with $\|x\|_{w,\infty}$. We write $\|A\|_{v,w,1}$ for the norm of
$A$ as an operator from  $d(v,1)$ into $d(w,1)$, and
$\|A\|_{w,v,\infty}$ for the norm of $A$ as an operator from
$e(w,\infty)$ into $e(v,\infty)$, and $\|A\|_{w,1}$ for the norm
of $A$ as an operator from $d(w,1)$ into itself, and
$\|A\|_{w,\infty}$ for the norm of $A$ as an operator from
$e(w,\infty)$ into itself.

\setcounter{theore}{0}
\begin{theor}[\!] Suppose $A=(a_{i,j})$ is a matrix operator
satisfying conditions $(1), (2)$ and $(3)$. If
\begin{equation*}
\sup\frac{S_n}{V_n}<\infty,
\end{equation*}
where $S_n=s_1+\cdots+s_n$ and $s_n=\sum_{k=1}^{\infty}w_ka_{k,n}$
and $V_n=v_1+\cdots+v_n${\rm ,} then $A$ is a bounded operator
from $d(v,1)$ into $d(w,1)${\rm ,} and also
\begin{equation*}
\|A\|_{v,w,1} =\sup_n\frac{S_n}{V_n}.
\end{equation*}
\end{theor}

\begin{proof}
By Proposition~2.1, it is sufficient to consider decreasing,
non-negative sequences. Let $x$ be in $d(v,1)$ such that $x_1\geq
x_2\geq\cdots\geq 0$ and $M=\sup\frac{S_n}{V_n}$. Then
\begin{align*}
\|Ax\|_{w,1}&=\sum_{n=1}^{\infty}w_n \left(\sum_{k=1}^{\infty}a_{n,k}x_k\right)\\[.4pc]
&=\sum_{n=1}^{\infty}s_nx_n\\[.4pc]
&=\sum_{n=1}^{\infty}S_n(x_n-x_{n+1})\\[.4pc]
&\leq M\sum_{n=1}^{\infty}V_n(x_n-x_{n+1}).
\end{align*}
Also, we have
\begin{equation*}
\|x\|_{v,1}=\sum_{n=1}^{\infty}V_n(x_n-x_{n+1}).
\end{equation*}
Therefore
\begin{equation*}
\|Ax\|_{w,1}\leq M\|x\|_{v,1},
\end{equation*}
and hence
\begin{equation*}
\|A\|_{v,w,1}\leq M.
\end{equation*}
To show that the constant $M$ is the best possible constant in the
above inequality, we take $x_1=x_2=\cdots=x_n=1$ and $x_k=0$ for
all $k\geq n+1$. Then
\begin{equation*}
\|x\|_{v,1}=V_n,\quad \|Ax\|_{w,1}=S_n.
\end{equation*}
Therefore
\begin{equation*}
\|A\|_{v,w,1}=M.
\end{equation*}
$\left.\right.$\vspace{-1.5pc}

\hfill $\Box$
\end{proof}
In the following statement we obtain  norm of general matrix
operator from $e(w,\infty)$ into $e(v,\infty)$.

\begin{theor}[\!]
Suppose $A=(a_{i,j})$ is a matrix operator satisfying conditions
$(1), (2)$ and $(3)$. If
\begin{equation*}
\sup\frac{Z_n}{V_n}<\infty,
\end{equation*}
where $Z_n=z_1+\cdots+z_n$ and
$z_n=\sum_{k=1}^{\infty}w_ka_{n,k}${\rm ,} then $A$ is a bounded
operator from $e(w,\infty)$ into $e(v,\infty)${\rm ,} and also
\begin{align*}
\|A\|_{w,v,\infty}&= \sup_n\frac{Z_n}{V_n}.
\end{align*}
\end{theor}

\begin{proof}
By Proposition~2.1, it is sufficient to consider decreasing,
non-negative sequences. Let $x$ be in $e(w,\infty)$ such that
$x_1\geq x_2\geq\cdots\geq 0$ and $\|x\|_{w,\infty}=1$. Then
\begin{equation*}
X_n\leq W_n,\quad \forall n.
\end{equation*}
Let $y=Ax$ and $c_{n,j}=\sum_{i=1}^{n}a_{i,j}$. We have
\begin{align*}
Y_n=\sum_{i=1}^{n}y_i&= \sum_{i=1}^{n}\sum_{j=1}^{\infty}a_{i,j}x_j\\[.4pc]
&= \sum_{j=1}^{\infty}c_{n,j}x_j\\[.4pc]
&= \sum_{j=1}^{\infty}(c_{n,j}-c_{n,j+1})X_j\\[.4pc]
&\leq \sum_{j=1}^{\infty}(c_{n,j}-c_{n,j+1})W_j\\[.4pc]
&= Z_n.
\end{align*}
If $C=\sup\frac{Z_n}{V_n},$ then
\begin{equation*}
\sup\frac{Y_n}{V_n}\leq C,
\end{equation*}
and hence $\|A\|_{w,v,\infty}\leq C.$

Since  $w\in e(w,\infty)$, $\|w\|_{w,\infty}=1$ and
$\|A(w)\|_{v,\infty}=C$, we have
\begin{equation*}
\|A\|_{w,v,\infty}=C.
\end{equation*}
If $A$ is a bounded matrix operator from $e(w,\infty)$ into
$e(v,\infty),$ then $A^t,$ the transpose matrix of $A,$ is a
bounded matrix operator of $d(v,1)$ into $d(w,1)$ and
\begin{equation*}
\|A^t\|_{v,w,1}=\|A\|_{w,v,\infty}.
\end{equation*}
$\left.\right.$\vspace{-1.5pc}

\hfill $\Box$
\end{proof}
Let $(a_n)$ be a non-negative sequence with $a_1>0$, and
$A_n=a_1+\cdots+a_n$. The N\"{o}rlund matrix $N_a=(a_{n,k})$ is
defined as follows:
\begin{align*}
a_{n,k}=\left\{\begin{array}{ll}
\frac{a_{n-k+1}}{A_n}, &1\le k\le n\\[.4pc]
0, &k>n
\end{array}.\right.
\end{align*}
If $\alpha\ge 0,$ the Cesaro matrix $C(\alpha)$ is matrix $N_a$
with
\begin{align*}
a_{n}=\left(\begin{array}{c} n+\alpha-2\\[.2pc] n-1
\end{array}\right).
\end{align*}
The Copson matrix of order $\alpha$ is the transpose matrix of
$C(\alpha)$, and we denote it with $C^t(\alpha)$. Also we denote
$C=C(1)$ and $C^t=C^t(1)$.

In the following statements, we consider the norm of Cesaro and
Copson matrices. It is enough to consider the sequence
$\big(\frac{s_n}{v_n}\big)$ instead of
$\big(\frac{S_n}{V_n}\big),$ because of the well-known fact listed
in the following lemma.

\setcounter{theore}{0}
\begin{lem}
If $m\leq\frac{s_n}{v_n}\leq M$ for all $n,$ then
$m\leq\frac{S_n}{V_n}\leq M$ for all $n.$
\end{lem}

\begin{proof}
It is elementary.\hfill $\Box$
\end{proof}

\setcounter{theore}{0}
\begin{propo}$\left.\right.$\vspace{.5pc}

\noindent If $w_n=\frac{1}{n}$ and $v_n=\frac{1}{n+\alpha}$ with
$\alpha\ge0,$ then $C(2)$ is a bounded operator from $d(v,1)$ into
$d(w,1)$ and also $C^t(2)$ is a bounded operator from
$e(w,\infty)$ into $e(v,\infty),$ and
\begin{equation*}
\|C(2)\|_{v,w,1}=\|C^t(2)\|_{w,v,\infty}=2(\alpha+1).
\end{equation*}
\end{propo}

\begin{proof}
We show that $\frac{s_n}{v_n}\le\frac{s_1}{v_1}$ for all $n$.
Therefore applying Lemma~3.1, we deduce that
$\frac{S_n}{V_n}\le\frac{S_1}{V_1}=s_1(\alpha+1),$ and by
Theorem~3.1, we have
\begin{equation*}
\|C(2)\|_{v,w,1}=2(\alpha+1).
\end{equation*}
Since
\begin{equation*}
s_1=\sum_{k=1}^{\infty}\frac{1}{\frac{1}{2}k(k+1)}=2,
\end{equation*}
for all $n$,
\begin{align*}
\frac{s_n}{v_n}&= (n+\alpha)\sum_{k=n}^{\infty}\frac{1}{\frac{1}{2}k(k+1)}\frac{k-n+1}{k}\\[.4pc]
&\le 2(n+\alpha)\sum_{k=n}^{\infty}\frac{1}{k(k+1)}\\[.4pc]
&\leq 2n\sum_{k=n}^{\infty}\frac{1}{k(k+1)}+2\alpha\sum_{k=1}^{\infty}\frac{1}{k(k+1)}\\[.4pc]
&= 2+2\alpha=\frac{s_1}{v_1}.
\end{align*}
This establishes the proof of the proposition.\hfill $\Box$
\end{proof}

\begin{propo}$\left.\right.$\vspace{.5pc}

\noindent If $w_n=\frac{1}{n}$ and $v_n=\frac{1}{n^{\alpha}}$ with
$0\le\alpha\leq1,$ then $C(2)$ is a bounded operator from $d(v,1)$
into $d(w,1)$ and also $C^t\!(2)$ is a bounded operator from
$e(w,\infty)$ into $e(v,\infty),$ and
\begin{equation*}
\|C(2)\|_{v,w,1}=\|C^t\!(2)\|_{w,v,\infty}=2.
\end{equation*}
\end{propo}

\begin{proof}
We show that $\frac{s_n}{v_n}\leq2$ for all $n$. Therefore
applying Lemma~3.1, we deduce that $\frac{S_n}{V_n}\leq2$. For all
$n$,
\begin{align*}
\frac{s_n}{v_n} &= n^{\alpha}\sum_{k=n}^{\infty}\frac{1}{\frac{1}{2}k(k+1)}\frac{k-n+1}{k}\\[.4pc]
&\le 2n^{\alpha}\sum_{k=n}^{\infty}\frac{1}{k(k+1)}\\[.4pc]
&\leq 2.
\end{align*}
Since $\frac{s_1}{v_1}=2,$ we have $\sup\frac{S_n}{V_n}=2.$ This
completes the proof of the proposition.\hfill $\Box$
\end{proof}

\setcounter{theore}{0}
\begin{coro}$\left.\right.$\vspace{.5pc}

If
\begin{equation*}
\sup_{n}\frac{1}{V_n}\sum_{k=1}^{n}\frac{W_k}{k}<\infty,
\end{equation*}
then the Cesaro matrix $C$ is a bounded operator from
$e(w,\infty)$ into $e(v,\infty),$ and
\begin{equation*}
\|C\|_{w,v,\infty}=\sup_{n}\frac{1}{V_n}\sum_{k=1}^{n}
\frac{W_k}{k}.
\end{equation*}
\end{coro}

\begin{proof}
By Theorem~3.1, we have
\begin{equation*}
\|C^t\|_{v,w,1}=\sup\frac{S_n}{V_n}.
\end{equation*}
Since $s_n=\frac{W_n}{n}$ and $\|C^t\|_{v,w,1}=
\|C\|_{w,v,\infty}$, we have the corollary.\hfill $\Box$
\end{proof}

\setcounter{theore}{2}
\begin{theor}[\!]
Suppose
\begin{equation*}
r=\sup\frac{W_n}{nv_n}<\infty.
\end{equation*}
Then the Copson operator $C^t$ maps $d(v,1)$ into $d(w,1)$ and
\begin{equation*}
\|C^t\|_{v,w,1}\leq r.
\end{equation*}
\end{theor}

\begin{proof}
Since $s_n=\frac{W_n}{n},$ we have $\sup\frac{S_n}{v_n}\leq r.$
Hence
\begin{equation*}
\|C^t\|_{v,w,1}=\sup\frac{S_n}{V_n}\leq r.
\end{equation*}
$\left.\right.$\vspace{-1.5pc}

\hfill $\Box$
\end{proof}

\begin{theor}[\!]
Suppose $v_n=\frac{1}{n^\alpha}$ and $W_n=n^{1-\alpha}$ with
$0\leq\alpha\le1.$ Then the Copson operator $C^t$ maps $d(v,1)$
into $d(w,1),$ and
\begin{equation*}
\|C^t\|_{v,w,1}=1.
\end{equation*}
Therefore
\begin{equation*}
\|C\|_{w,v,\infty}=1.
\end{equation*}
\end{theor}

\begin{proof}
For all $n$, $\frac{W_n}{nv_n}=1$, and  therefore $r=1$. Hence
\begin{equation*}
\|C^t\|_{v,w,1}\leq 1.
\end{equation*}
Since $\frac{s_1}{v_1}=1$, we deduce that
\begin{equation*}
\|C^t\|_{v,w,1}=1.
\end{equation*}
$\left.\right.$\vspace{-1.5pc}

\hfill $\Box$
\end{proof}

\begin{theor}[\!]
Suppose $w_n=\frac{1}{n^\alpha}$ and $V_n=n^{1-\alpha}$ with
$0\leq\alpha\leq1.$ Then the Cesaro operator C maps $d(v,1)$ into
$d(w,1),$ and
\begin{equation*}
\|C\|_{v,w,1}\leq\frac{1}{1-\alpha}\zeta(1+\alpha).
\end{equation*}
Therefore the Copson operator $C$ maps $e(w,\infty)$ into
$e(v,\infty),$ and
\begin{equation*}
\|C^t\|_{w,v,\infty}\leq\frac{1}{1-\alpha}\zeta(1+\alpha).
\end{equation*}
\end{theor}

\begin{proof}
By mean value theorem for all $n$, we have
\begin{equation*}
\frac{1-\alpha}{n^{\alpha}}\leq n^{1-\alpha}-(n-1)^{1-{\alpha}}.
\end{equation*}
Since $v_n=n^{1-\alpha}-(n-1)^{1-\alpha},$
\begin{equation*}
\frac{s_n}{v_n}\leq\frac{n^\alpha}{1-\alpha}s_n,
\end{equation*}
and hence $\sup\frac{s_n}{v_n}\leq\frac{1}{1-\alpha}\sup n^\alpha
s_n.$

The sequence $(n^\alpha s_n)$ is decreasing (Lemma~2.7 of
\cite{6}), and therefore
\begin{equation*}
\sup\frac{s_n}{v_n}\leq\frac{1}{1-\alpha}s_1=\frac{1}{1-\alpha}
\zeta(1+\alpha).
\end{equation*}
This completes the proof of the theorem.\hfill $\Box$
\end{proof}

We recall that the Hilbert operator $H$ is defined by the matrix
\begin{equation*}
a_{i,j}=\frac{1}{i+j}.
\end{equation*}

\setcounter{theore}{1}
\begin{lem}
If $0\leq\alpha\leq1,$ then
\begin{equation*}
\sup_n{n^\alpha}\sum_{k=1}^{\infty}\frac{1}{k^{\alpha}(k+n)}
=\frac{\pi}{\sin\alpha\pi}.
\end{equation*}
\end{lem}

\begin{proof}
It is elementary.\hfill $\Box$
\end{proof}
In the following statement, we consider the upper bound of $H$.

\setcounter{theore}{5}
\begin{theor}[\!]
Suppose $w_n=\frac{1}{n^\alpha}$ and $V_n=n^{1-\alpha}$ where
$0\leq\alpha\leq1.$ Then the Hilbert matrix operator $H$ maps
$d(v,1)$ into $d(w,1),$ and
\begin{equation*}
\|H\|_{v,w,1}\leq\frac{\pi}{(1-\alpha)\sin\alpha\pi}.
\end{equation*}
Therefore the Hilbert operator $H$ maps $e(w,\infty)$ into
$e(v,\infty),$ and
\begin{equation*}
\|H\|_{w,v,\infty}\leq\frac{\pi}{(1-\alpha)\sin\alpha\pi}.
\end{equation*}
\end{theor}

\begin{proof} We have
$\frac{s_n}{v_n}\leq\frac{n^\alpha}{1-\alpha}s_n$ which is similar
to the previous theorem. Applying Lemma~3.2, we have
$\sup{n^\alpha}s_n=\frac{\pi}{\sin\alpha\pi}$, and so
\begin{equation*}
\frac{s_n}{v_n}\leq{\frac{\pi}{(1-\alpha)\sin\alpha\pi}}.
\end{equation*}
This completes the proof of the theorem.\hfill $\Box$
\end{proof}

\end{document}